\documentclass[12pt,doublespace]{amsart}
\usepackage{amssymb}
\newtoks\prt
\newtheorem{thm}{Theorem}
\newtheorem{lemma}[thm]{Lemma}

\theoremstyle{definition}

\newtheorem{remark}[thm]{Remark}

\def\ep{\varepsilon}

\def\H{\mathcal H}

\def\diam{\operatorname{diam}}

\def\md{\operatorname{md}}

\newcommand{\R}{{\mathbb{R}}}

\usepackage{amssymb}
\usepackage{amsmath}
\usepackage{amsthm}

\newcommand{\e}{\varepsilon}

\newtoks\by
\newtoks\paper
\newtoks\book
\newtoks\jour
\newtoks\yr
\newtoks\pages
\newtoks\vol
\newtoks\publ
\newtoks\eds
\newtoks\proc
\newtoks\mathrev
\def\ota{{\hbox{???}}}
\def\cLear{\by=\ota\paper=\ota\book=\ota\jour=\ota\yr=\ota
\pages=\ota\vol=\ota\publ=\ota}
\def\endpaper{\the\by, \textit{\the\paper},
{\the\jour} \textbf{\the\vol} (\the\yr), \the\pages.\cLear}
\def\endbook{\the\by, \textit{\the\book}, \the\publ, \the\yr.\cLear}
\def\endprep{\the\by, \textit{\the\paper}, \the\jour.\cLear}
\def\endproc{\the\by, \textit{\the\paper}, \the\book,
\the\publ, \the\yr, \the\pages.\cLear}

\begin{document}

\title
{Morse-Sard theorem for d.\,c. curves}

\author{D. Pavlica}
\address{Mathematical Institute, Academy of Science of the Czech Republic\\ \v Zitn\'a 25\\
115 67 Praha 1\\
Czech Republic}

\email{pavlica@karlin.mff.cuni.cz}


\begin{abstract}
Let $f:I\to X$ be a d.c. mapping, where $I\subset \R$ is an open interval
 and $X$ a Banach space. Let $C_f$ be the set of critical points of $f$.
 We prove that $f(C_f)$ has zero $1/2$-dimensional Hausdorff measure.

\end{abstract}



\maketitle


Let $Z$ and $X$ be Banach spaces, $U\subset Z$ be an open convex
set and $f\colon U \to X$ be a mapping. We say that $f$ is a {\em
d.c.\  mapping\/} if there exists a continuous convex function $h$
on $U$ such that $y^* \circ f + h$ is a continuous convex function
for each $y^* \in Y^*, \|y^*\|=1$. We say that $f:U \to X$ is locally d.c.
if for each $x\in U$ there exists open convex $U'$ such that $x\in U'\subset U$ and 
$f|_{U'}$ is d.c.

This notion of d.c. mappings between Banach spaces (see \cite{VZ}) generalizes Hartman's \cite{H} notion
of  d.c. mapping between Euclidean spaces.
 Note that in this case 
 it is easy to see that $F$ is d.c. if and only if all its components
  are d.c. (i.e. are differences of two convex functions).

\smallskip

For $f:U\to X$ we denote
$C_f:=\{x\in U: f'(x)=0\}$.

A special case of  \cite[Theorem 3.4.3]{F}  says that for a
mapping $f:\R^m \to X$ of class $C^2$, where $X$ is a normed vector
space, the set $f(C_f)$ has zero $(m/2)$-dimensional
Hausdorff measure. 

We will generalize this result in case $m=1$ showing that it
is sufficient to
 suppose that $f$ is d.c. on $I$ (equivalently: $f$ is continuous and
 $f'_+$ is locally of bounded
 variation on $I$).

Similar generalization of above mentioned result on $C^2$ mappings holds for $m=2$ as is shown (by completely different method)
in \cite{PZ} where it is proved that $f(C_f)$ has zero 1-dimensional Hausdorff measure for d.c. mapping $f:\R^2 \to X$.

Whether $f(C_f)$ has zero $(m/2)$-dimensional Hausdorff measure for each d.c. mapping $f:\R^m \to X$ for $m\ge 3$
remains open even for $X$ an Euclidean space. 

\smallskip 

We denote  $\alpha$-dimensional Hausdorff measure (on a metric space $X$) by $\H^{\alpha}$
 and for each $Y \subset X$ we put (see \cite{M})
$$ \H^{\alpha}_{\infty}(Y) = \frac{\omega_{\alpha}}
{2^{\alpha}}\cdot   \inf \{ \sum_{i=1}^{\infty} \diam^{\alpha}
(M_i):\ Y \subset
 \bigcup_{i=1}^{\infty} M_i\},$$
 where $\omega_{\alpha}= (\Gamma(1/2))^{\alpha}\cdot (\Gamma(\alpha/2 +
 1))^{-1}$.

For an open interval $I$, a Banach space $X$,  $g:I\to X$
 and $x\in I$, we denote
$$\md(g,x):=\lim_{r\to 0}\frac{\|g(x+r)-g(x)\|}{|r|},\quad
x\in I.$$

If $g$ is Lipschitz, then $\md(g,x)$ exists a.e. on $I$. This fact is a
special case of Kirchheim's theorem \cite[Theorem 2]{Ki} on a.e.
 metric differentiability of Lipschitz mappings (from $\R^n$ to $X$).
In a standard way, we obtain the following more general fact.

\begin{lemma}
Let $I$ be an open interval, $X$ be a Banach space $X$, and let $g:I\to X$
have bounded variation on $I$. Then $\md(g,x)$ exists
almost everywhere on $I$.
\end{lemma}
\begin{proof}
We may suppose $I=\R$. Denote $s(x) = \bigvee_{0}^x g,\ x \in \R$.
By \cite[2.5.16.]{F} there exists a Lipschitz mapping $H: \R \to
X$
 such that $g = H \circ s$. By \cite[Theorem 2]{Ki}, $\md(H,x)$ exists
 a.e. on $\R$. Now, changing in the obvious way the last argument
 of \cite[2.9.22.]{F}, we obtain our assertion.
\end{proof}

\begin{thm}\label{1/2}
Let $X$ be a Banach space, $I\subset \R$  be an open interval and
$f: I \to X$ a locally  d.c. mapping. Let $C := \{x \in I:\
f'(x) = 0\}$. Then $\H^{1/2}(f(C))=0$.
\end{thm}
\begin{proof}
Note that $f$ is continuous on $I$  (see \cite[Proposition 1.10.]{VZ}).
 By  \cite[Theorem 2.3]{VZ}, $f'_+$ exists
and has locally bounded variation on $I$. Consider an arbitrary
 interval $[a,b] \subset I$. It is clearly sufficient to
 prove $\H^{1/2}(f(C_1))=0$, where $C_1 := C \cap (a,b)$.

Let $N_1$ be the set of all isolated points of $C_1$ and
 $N_2:=\{x\in C_1:\  \md(f'_+,x)\  \text{does not exist} \}$.
Since $N_1$ is countable, $\H^{1/2}(f(N_1))=0$.

To prove $\H^{1/2}(f(N_2))=0$, consider an arbitrary $\ep >0$.
 By Lemma 1, we find  a countable disjoint
 system of open intervals $\{(a_i,b_i):\ i \in J\}$ such that
$$ N_2 \subset \bigcup_{i \in J} (a_i,b_i) \subset (a,b),\
 \sum_{i\in J} (b_i-a_i) < \ep\ \ \text{and}\ \ (a_i,b_i) \cap N_2 \neq
 \emptyset,\ i \in J.$$
Clearly  $\|f'_+(x)\| \leq \bigvee_{a_i}^{b_i} \, f'_+$ for each $i \in J$
 and $x \in (a_i,b_i)$. Using continuity of $f$ and  \cite[Chap. I, par. 2,
 Proposition 3]{B}, we obtain
$$   \diam(f((a_i,b_i))) \leq  (b_i - a_i) \cdot  \bigvee_{a_i}^{b_i}
\, f'_+.$$
 Therefore, using the  Cauchy-Schwartz inequality, we obtain
\begin{eqnarray*}
\H^{1/2}_\infty(f(N_2))&\le&
 \frac{\omega_{1/2}}{2^{1/2}} \sum_{i\in J}
((b_i - a_i) \cdot  \bigvee_{a_i}^{b_i}
\, f'_+)^{1/2}\\
 &\le& \omega_{1/2} (\sum_{i\in J} (b_i-a_i))^{1/2} (\sum_{i\in J}
\bigvee_{a_i}^{b_i} f'_+)^{1/2} \le
\omega_{1/2}\,  \ep^{1/2} \left(\bigvee_a^b
f'_+\right)^{1/2}  .
 \end{eqnarray*}
 Since $\e>0$ is arbitrary, we have $\H^{1/2}_\infty(f(N_2))=0$,
 consequently
 (see \cite[Lemma 4.6.]{M})  we obtain
$\H^{1/2}(f(N_2))=0$.

To finish the proof, it is sufficient to prove
$\H^{1/2}(f(C_2))=0$, where $C_2 = C_1 \setminus (N_1 \cup N_2)$.
Let $\e >0$ be arbitrary. Clearly $\md(f'_+,x)=0$ for each $x \in C_2$.
 Therefore, for
each $x\in C_2$, we can choose $\delta_x>0$ such that
 $[x-\delta_x,x+ \delta_x] \subset (a,b)$ and
 $\|f'_+(y)\| \leq
 \ep\, |y-x|$ for each $y \in [x-\delta_x,x+ \delta_x]$.
 Using continuity of $f$ and \cite[Chap. I, par. 2, Proposition 3]{B}, we
 obtain  \ $\diam(f([x-\delta_x,x+\delta_x])) \leq  2 \ep (\delta_x)^2$.

 Besicovitch's Covering Theorem (see \cite{M}) easily implies that
 we can choose a countable set $A \subset C_2$ such that
$$ C_2 \subset \bigcup_{x \in A} [x-\delta_x, x +\delta_x]\ \ \text{and}\ \
 \sum_{x\in A} 2 \delta_x \leq c\, (b-a),$$
 where  $c$ is an absolute constant (not depending on $\ep$).
Since
$\e>0$ is arbitrary,
$$f(C_2)\subset \bigcup_{x\in A} f([x-\delta_x,x+\delta_x]),$$
 and
$$
\sum_{x\in A} (\diam(f([x-\delta_x,x+\delta_x])))^{1/2} \leq
 \sum_{x\in A} \sqrt{2 \ep}\, \delta_x \leq \sqrt{2 \ep}\, c (b-a),$$
 we have $\H^{1/2}_\infty(f(C_2))=0$, hence
$\H^{1/2}(f(C_2))=0$.
             \end{proof}

\begin{remark}
Since each $C^2$-function on $I$ is a locally d.c. function (see \cite{VZ}),  
 \cite[3.4.4.]{F} implies that the conclusion of
 Theorem \ref{1/2} does not hold with $\H^{\alpha}$ $(\alpha < 1/2)$ in general. 
\end{remark}

\end{document}